\documentclass[10pt]{amsart}
\usepackage{amsthm,amsmath,amssymb,amsfonts}
\usepackage[american]{babel}
\usepackage{cite}
\usepackage{url}

\newtheorem{theorem}{Theorem}
\newtheorem*{theorem*}{Theorem}

\newtheorem{proposition}[theorem]{Proposition}

\theoremstyle{definition}

\theoremstyle{remark}

\begin{document}


\newcommand{\NU}[1]{\mbox{\rm V}(#1)}
\newcommand{\lara}[1]{\langle{#1}\rangle}
\newcommand{\ZZ}{\mbox{$\mathbb{Z}$}}
\newcommand{\sZZ}{\mbox{$\scriptstyle\mathbb{Z}$}}
\newcommand{\QQ}{\mbox{$\mathbb{Q}$}}
\newcommand{\paug}[2]{\varepsilon_{#1}(#2)}


\title[Orders of torsion units in integral group rings]
{The orders of torsion units in integral group rings
of finite solvable groups}

\author{Martin Hertweck}
\address{Universit\"at Stuttgart, Fachbereich Mathematik, IGT,
Pfaffenwaldring 57, 70550 Stuttgart, Germany}
\email{hertweck@mathematik.uni-stuttgart.de}



\date{\today}


\begin{abstract}
It is shown that for any torsion unit of augmentation one in the 
integral group ring $\ZZ G$ of a finite solvable group $G$,
there is an element of $G$ of the same order.
\end{abstract}

\maketitle

\section{Introduction}\label{Sec:intro}

Does a torsion unit in the integral group ring $\ZZ G$ of a 
finite group $G$ necessarily has the same order as some element of $G$?
Of course, torsion coming from the coefficient ring should be excluded,
that is, only torsion units from $\NU{\ZZ G}$, the group of units of 
augmentation one in $\ZZ G$, will be considered. When $G$ is solvable,  
this note provides an affirmative answer, but else we are left without 
a clue. It is known that the answer is affirmative for torsion units of
prime power order, and that the order of any torsion unit divides the 
exponent of $G$, see Cohn and Livingstone (1965) or Zassenhaus (1974).
The question to be considered is included as Research Problem~8 in 
the book of Sehgal (1993), where an affirmative answer
is given for metabelian groups by Lemma (37.14).
It should be seen in the context of the stronger 
conjectures of Zassenhaus about rational conjugacy of torsion units
(see Sehgal, 1993, Chapter~5).
In the course of his work on these conjectures, Weiss (1988) proved a 
beautiful theorem about $p$-permutation lattices 
(cf.\ Sehgal, 1993, Appendix).
The present contribution relies on this theorem, as does already 
previous work of the author (Hertweck, 2006, 2007)
from which the method of proof is taken.

Recall that for a group ring element $u=\sum_{g\in G}a_{g}g$ 
(all $a_{g}$ in $\ZZ$), its partial augmentation
with respect to an element $x$ of $G$, or rather its
conjugacy class $x^{G}$ in $G$, is the sum $\sum_{g\in x^{G}}a_{g}$;
we will denote it by $\paug{x}{u}$.
Our result is as follows.
\begin{theorem*}
Let $G$ be a finite solvable group. Then any torsion unit in $\NU{\ZZ G}$
has a nonzero partial augmentation with respect to a conjugacy class of a
group element of the same order. In particular, the orders of the torsion 
units in $\NU{\ZZ G}$ are the orders of the elements of $G$.
\end{theorem*}

The author's interest in proving such a result was stimulated by 
recent work of H\"ofert and Kimmerle on the prime
graph of $\NU{\ZZ G}$ (see Kimmerle, 2006, \S4).
They showed, for solvable $G$, that a unit in $\NU{\ZZ G}$ of order
$pq$, for distinct primes $p$ and $q$, can only exist if $G$ has an
element of order $pq$. Their method of proof is the one described in
Cohn and Livingstone (1965) and Zassenhaus (1974),
which provides certain congruences for the partial augmentations of the 
$p$-th and $q$-th powers of the unit in question.
Roughly speaking, they noticed that in a specific situation, which
can be reached inductively, the presence of a normal $p$-subgroup
of $G$ can be used to turn one of these congruences into an equality, 
yielding the result.

Finally, we remark that, as with other results in this field, the 
theorem can be formulated for more general coefficient rings than
$\ZZ$, notably for the semilocalization of $\ZZ$ at the prime divisors
of the order of $G$.
We stick to the basic case only to reduce notational effort.

\section{Proof of the theorem}\label{Sec:proof}

The information provided by Weiss' theorem which we will use 
can be found in \S4 of Hertweck (2006) and is described in the following 
proposition. Let $G$ be a finite group. We fix a rational prime $p$
and let $\ZZ_{p}$ stand for the $p$-adic integers.

\begin{proposition}\label{P1}
Suppose that $G$ has a normal $p$-subgroup $N$, and that $U$ is a finite
subgroup of $\NU{\ZZ G}$ which maps to $1$ under the natural map
$\ZZ G\rightarrow \ZZ G/N$. Then $U$ is a $p$-group, and
$\ZZ_{p}G$, when considered as $\ZZ_{p}U$-lattice via the (right)
multiplication action of $U$, is projective.
\end{proposition}

Indeed, under the assumptions made, it readily follows that $U$ is a 
$p$-group, and Weiss' theorem guarantees that $\ZZ_{p}G$ is a permutation
lattice for $U$ over $\ZZ_{p}$. From that, one obtains that $U$ is 
conjugate to a subgroup of $N$ by a unit in $\QQ G$ (see 
Hertweck, 2006, Proposition~4.2).
Combined together, this implies the proposition.

The next proposition describes how this will be applied.
The line of proof is that of Hertweck (2007), Proposition~2.2,
with roles of torsion unit and group element interchanged,
but we repeat the argument for convenience of the reader.

\begin{proposition}\label{P2}
Suppose that $G$ has a normal $p$-subgroup $N$, and that $u$ is a 
torsion unit in $\NU{\ZZ G}$ whose image under the natural map
$\ZZ G\rightarrow \ZZ G/N$ has strictly smaller order than $u$.
Then $\paug{g}{u}=0$ for every element $g$ of $G$ whose 
$p$-part has strictly smaller order than the $p$-part of $u$. 
\end{proposition}
\begin{proof}
By assumption, $\lara{u}$ has a subgroup of order $p$ which maps to $1$
under the map $\ZZ G\rightarrow \ZZ G/N$. Suppose that $g$ is an element
of $G$ whose $p$-part has strictly smaller order than the $p$-part of $u$. 
Let $C$ be an (abstract) cyclic group whose order is the least common 
multiple of the orders of $u$ and $g$.
Set $M=\ZZ_{p}G$, viewed as $\ZZ_{p}C$-lattice by letting a generator $c$ 
of $C$ act by $m\cdot c=g^{-1}mu$ for $m\in M$. By (38.12) in Sehgal (1993),
we have to show that $\chi(c)=0$ for the character $\chi$ of $C$ afforded
by $M$. Note that $c$ is $p$-singular, so this will follow once we have
shown the stronger statement that $M$ is a projective $\ZZ_{p}C$-lattice, 
by Green's Theorem on Zeros of Characters (see (19.27)
in Curtis and Reiner, 1981).
This in turn follows from the assumption, meaning that 
the action of the subgroup $P$ of order $p$ in $C$ is given
by a multiplication action of the subgroup of order $p$ in $\lara{u}$,
which shows that $M$ is a projective 
$\ZZ_{p}P$-lattice, by Proposition~\ref{P1}.
We provide details. Set $k=\ZZ_{p}/p\ZZ_{p}=\ZZ/p\ZZ$.
It is enough to show that $k\otimes_{\sZZ_{p}}M$ as $kC$-module is projective
(see (20.10) or (30.11) in Curtis and Reiner, 1981), i.e., that
$k\otimes_{\sZZ_{p}}M$ is projective relative to a Sylow $p$-subgroup of $C$.
It is well-known that this follows from the projectivity of
$k\otimes_{\sZZ_{p}}M$ as $kP$-module (see Hertweck, 2006, Lemma~3.2).
\end{proof}

We will use the following elementary observation. 
Let $N$ be a normal subgroup of $G$ and set $\bar{G}=G/N$. 
For a torsion unit $u$ in $\ZZ G$, we shall extend the bar convention 
when writing $\bar{u}$ for the image of $u$ under the natural map 
$\ZZ G\rightarrow \ZZ\bar{G}$. 
We let $\sim$ denote the conjugacy relation in a group.
Since any conjugacy class of $G$ maps
onto a conjugacy class of $\bar{G}$, we have for any $x\in G$: 
\begin{equation}\label{E1}
\paug{\bar{x}}{\bar{u}} = \sum_{g^{G}:\;\bar{g}\sim\bar{x}}
\paug{g}{u}. 
\end{equation}

We now turn to the proof of the theorem. Let $G$ be a finite solvable
group, and let $u$ be a torsion unit in $\NU{\ZZ G}$.
The proof is by induction on the order of $G$. Since $G$ is solvable,
we can choose a nontrivial normal $p$-subgroup $N$ of $G$. Set
$\bar{G}=G/N$. Then $\bar{u}$ denotes the image of $u$ under the 
homomorphism $\NU{\ZZ G}\rightarrow\NU{\ZZ\bar{G}}$. Inductively, there
is an element $x$ in $G$ such that $\bar{x}$ and $\bar{u}$ have the 
same order and $\paug{\bar{x}}{\bar{u}}\neq 0$. Let $g\in G$ such that
$\bar{g}\sim\bar{x}$. Note that if $g$ has strictly greater order than
$u$, then $\paug{g}{u}=0$ by Theorem~2.3 in Hertweck (2007).
We distinguish two cases. Firstly, suppose that $\bar{u}$ has the same
order as $u$. Then the order of $g$ is strictly greater than or equal 
to the order of $u$. Thus the sum in \eqref{E1} extends only over those
classes $g^{G}$ for which $g$ has the same order as $u$, and since at 
least one summand must contribute to its nonzero value, we obtain the 
desired result. Secondly, suppose that $\bar{u}$ has strictly smaller 
order than $u$. Since $\bar{g}$ and $\bar{u}$ have the same order, the
orders of $g$ and $u$ have the same $p^{\prime}$-part. Hence if $g$ has 
strictly smaller order than $u$, then $\paug{g}{u}=0$ by 
Proposition~\ref{P2}. Thus the sum in \eqref{E1} again extends only over
those classes $g^{G}$ for which $g$ has the same order as $u$, 
and we are done.

\section*{References}

\begin{list}{}{
\setlength{\leftmargin}{\parindent}
\setlength{\labelwidth}{\parindent}
\addtolength{\labelwidth}{-\labelsep}}
\item[Cohn, J.~A., Livingstone, D. (1965).]
On the structure of group algebras.~I.
{\it Canad. J. Math.} 17:583--593.
\item[Curtis, C.~W., Reiner, I. (1981).]
{\it Methods of representation theory.
With applications to finite groups and orders. Vol.~I}.
New York: John Wiley \& Sons Inc.
\item[Hertweck, M. (2006).] 
On the torsion units of some integral group rings.
{\it Algebra Colloq.} 13:329--348.
\item[Hertweck, M. (2007).] 
Partial augmentations and Brauer character values of
torsion units in group rings.
{\it Comm.\ Algebra}, to appear (e-print \url{arXiv:math.RA/0612429v2}).
\item[Kimmerle, W. (2006).]
On the prime graph of the unit group of integral group rings of finite groups.
In: Chin, W., Osterburg, J., Quinn, D., eds. 
{\em Groups, Rings and Algebras,
A conference in Honor of Donald S.\ Passman. Contemp. Math. Vol.~420}.
Providence, RI: Amer. Math. Soc. pp.~215--228.
\item[Sehgal, S.~K. (1993).]
{\it Units in integral group rings}.
Pitman Monographs and Surveys in Pure and Applied Mathematics Vol.~69.
Harlow: Longman Scientific \& Technical.
\item[Weiss, A. (1988).]
Rigidity of {$p$}-adic {$p$}-torsion.
{\it Ann. of Math. (2)} 127:317--332.
\item[Zassenhaus, H. (1974).]
On the torsion units of finite group rings.
In: {\it Studies in mathematics (in honor of A. Almeida Costa).} 
Lisbon: Instituto de Alta Cultura, pp.~119--126.
\end{list}

\end{document}